\documentclass{article}
\pdfoutput=1

\usepackage{amsmath}
\usepackage{amsfonts}
\usepackage{graphicx}
\usepackage{fancyhdr}

\usepackage[colorlinks=true,%
citecolor=black,%
filecolor=black,%
linkcolor=black,%
urlcolor=blue,%
pdfauthor={S. Congreve, J. Gedicke, and I. Perugia},%
pdftitle={Numerical investigation of the conditioning for plane wave discontinuous Galerkin methods}]{hyperref}

\title{Numerical investigation of the conditioning for plane wave discontinuous Galerkin methods\thanks{\textbf{accepted in Numerical Mathematics and Advanced Applications - ENUMATH 2017, Lecture Notes in Computational Science and Engineering, Vol. 126, Springer, 2018}}}

\author{
Scott Congreve\thanks{University of Vienna, Faculty of Mathematics, Oskar-Morgenstern-Platz 1, 1090 Vienna, Austria, {\tt scott.congreve@univie.ac.at}} 
\and 
Joscha Gedicke\thanks{University of Vienna, Faculty of Mathematics, Oskar-Morgenstern-Platz 1, 1090 Vienna, Austria, {\tt joscha.gedicke@univie.ac.at}}
\and 
Ilaria Perugia\thanks{University of Vienna, Faculty of Mathematics, Oskar-Morgenstern-Platz 1, 1090 Vienna, Austria, {\tt ilaria.perugia@univie.ac.at}}
}

\date{}

\renewcommand{\vec}[1]{\mathbf{#1}}

\begin{document}

\maketitle

\begin{abstract}
We present a numerical study to investigate the conditioning 
of the plane wave discontinuous Galerkin discretization of the Helmholtz problem.
We provide empirical evidence that the spectral condition number of
the plane wave basis on a single element
depends algebraically on the mesh size and the wave number, and
exponentially on the number of plane wave directions; we also test
its dependence on the element shape.
We show that the conditioning of the global system can be improved by orthogonalization of the
local basis functions with the modified Gram-Schmidt algorithm, which results in significantly fewer GMRES iterations for solving the discrete problem iteratively.
\end{abstract}

\section{Introduction}
\label{perugia_mini4:sec:1}
For the numerical approximation of the Helmholtz problem, 
it has been shown that, by using
non-polynomial basis functions, it is possible 
to reduce the pollution effect in finite element approximations.
One special class of such methods are Trefftz finite element methods, which use basis functions that are
local solutions of the homogeneous problem under consideration. 
For the Helmholtz problem, a common choice of Trefftz basis
functions are plane waves; when used in connection to a
discontinuous Galerkin variational framework, they lead to the plane
wave discontinuous Galerkin (PWDG) method 
\cite{PerugiaMS4_BM08,PerugiaMS4_CD98,PerugiaMS4_CD03,PerugiaMS4_GHP09,PerugiaMS4_HMP11,PerugiaMS4_HMP16b}.
Containing information on the oscillatory behaviour of the
solutions already in the approximation spaces, PWDG deliver better
accuracy than standard polynomial finite element methods, for a
comparable number of degrees of freedom.
In addition, they involve only evaluation of
basis functions on mesh interelement boundaries; 
hence, they can be easily used in connection with general polytopal
meshes.
However, it is well known that these basis functions are ill conditioned for 
small mesh sizes, small wavenumbers and large numbers of plane wave
directions \cite{PerugiaMS4_HMK02,PerugiaMS4_LHM13}. 

The aim of this paper is to numerically investigate the
dependence of the 
elemental and global
condition numbers of the PWDG
system matrix on the size and shape of the local (convex)
polygonal element,
the wavenumber,
and the number of plane wave directions in the local
approximation spaces.

\section{The PWDG method for the Helmholtz problem}
Let $\Omega\subset \mathbb{R}^2$ be a bounded Lipschitz domain and $k>0$ denote the wave\-number.
We consider the homogeneous Helmholtz problem with impedance boundary condition:
\begin{align}\label{perugia_problem}
\begin{split}
	-\Delta u -k^2 u &=  0 \quad \textrm{in } \Omega,\\
	\nabla u \cdot \vec{n} + iku &= g \quad \textrm{on } \partial\Omega,
\end{split}
\end{align}
where $i$ denotes the imaginary unit, $\vec{n}$ is the unit outward normal and $g\in L^2(\partial\Omega)$
is given.
The variational formulation of problem \eqref{perugia_problem} reads as follows:
find $u\in H^1(\Omega)$ such that 
\begin{align}\label{perugia_variational}
\int_\Omega (\nabla u\cdot \nabla \bar{v} - k^2 u \bar{v})dx + ik\int_{\partial\Omega} u\bar{v}\, \mathrm{d}s = \int_{\partial\Omega} g\bar{v}\, \mathrm{d}s
\quad \mbox{for all } v\in H^1(\Omega).
\end{align}
Problem \eqref{perugia_variational} is well posed by the Fredholm
alternative argument \cite{PerugiaMS4_Melenk}.
\par

We consider a shape-regular, uniform partition $\mathcal{T}_h$ of the domain $\Omega$ into convex
polygons $K\in\mathcal{T}_h$ of diameter $h$. We define the mesh skeleton $\mathcal{F}_h = \cup_{K\in\mathcal{T}_h}\partial K$, and denote the interior mesh skeleton by $\mathcal{F}_h^I = \mathcal{F}_h\backslash\partial\Omega$. 
For an element $K\in\mathcal{T}_h$ we define the plane wave space $\mathrm{PW}_p(K)$ of degree $p$ as
\[
  \mathrm{PW}_p(K) = \{ v\in L^2(K) \;:\; v(x) = \sum_{j=1}^p\alpha_j \exp(ik\vec{d}_j\cdot (\vec{x}-\vec{x}_K)), \alpha_j\in\mathbb{C} \},
\]
where $\vec{x}_K$ is the mass center of $K$, and $\vec{d}_j$, $|\vec{d}_j|=1$, $j=1,\ldots,p$, are $p$ unique directions.
Since, in general, small angles between those directions lead to bad conditioning of the basis,
we consider equally spaced directions. The 
PWDG space is defined as
\begin{align*}
  \mathrm{PW}_p(\mathcal{T}_h) = \{ v_{hp}\in L^2(\Omega) \;:\; v_{hp}|_K\in\mathrm{PW}_p(K)\quad\mbox{for all } K\in\mathcal{T}_h\}.
\end{align*} 
The functions in $\mathrm{PW}_p(\mathcal{T}_h)$ are local solutions of the homogeneous Helmholtz problem;
therefore, they exhibit the Trefftz property
\begin{align*}
  -\Delta v_{hp} - k^2 v_{hp} = 0\quad\mbox{for all } v_{hp}\in \mathrm{PW}_p(K).
\end{align*}
We assume uniform local resolution, i.e., we employ the same uniformly distributed directions $\vec{d}_j$, $j=1,\ldots,p$ on each element $K\in\mathcal{T}_h$.
\par

We use the standard notation for averages and
normal jumps of traces across inter-element
boundaries, namely $\{\!\!\{\cdot\}\!\!\}$ and $[\![\cdot]\!]$,
respectively, and 
denote by $\nabla_h$ the elementwise application of
$\nabla$.
Hence, we can formulate the PWDG method as follows \cite{PerugiaMS4_GHP09,PerugiaMS4_HMP11,PerugiaMS4_HMP16b}:
find $u_{hp}\in \mathrm{PW}_p(\mathcal{T}_h)$ such that
\begin{align}\label{perugia_pwdg}
\mathcal{A}_h(u_{hp},v_{hp})  = \ell_h(v_{hp})
\quad\mbox{for all } v_{hp} \in \mathrm{PW}_p(\mathcal{T}_h),
\end{align}
where 
\begin{align*}
\mathcal{A}_h(u_{hp},v_{hp}) &:= 
i\left[ -\int_{\mathcal{F}_h^I} \{\!\!\{u\}\!\!\} [\![\nabla_h\bar{v}]\!]\, \mathrm{d}s 
+ \int_{\mathcal{F}_h^I}\{\!\!\{\nabla_h u\}\!\!\}\cdot [\![\bar{v}]\!]\, \mathrm{d}s\right.\\
&\qquad-\frac{1}{2}\int_{\partial\Omega}u \nabla_h\bar{v}\cdot\vec{n} \, \mathrm{d}s
+\left.\frac{1}{2}\int_{\partial\Omega}\nabla_h u\cdot\vec{n}\bar{v}\, \mathrm{d}s\right]\\
&\qquad+\frac{1}{2k}\int_{\mathcal{F}_h^I}[\![\nabla_hu]\!][\![\nabla_h\bar{v}]\!]\, \mathrm{d}s
+\frac{k}{2}\int_{\mathcal{F}_h^I} [\![u]\!]\cdot[\![\bar{v}]\!]\, \mathrm{d}s\\
&\qquad+\frac{1}{2k}\int_{\partial\Omega}(\nabla_hu\cdot\vec{n})(\nabla_h\bar{v}\cdot\vec{n})\, \mathrm{d}s
+\frac{k}{2}\int_{\partial\Omega}u\bar{v}\, \mathrm{d}s,\\
\ell_h(v) &:= \frac{1}{2k}\int_{\partial\Omega}g\nabla_h\bar{v}\cdot\vec{n} \, \mathrm{d}s - \frac{i}{2}\int_{\partial\Omega}g\bar{v}\, \mathrm{d}s.
\end{align*}

The PWDG method \eqref{perugia_pwdg} is unconditionally well-posed and stable \cite{PerugiaMS4_BM08,PerugiaMS4_CD98}.
The $h$, $p$ and $hp$ convergence has been studied in \cite{PerugiaMS4_BM08,PerugiaMS4_GHP09,PerugiaMS4_HMP11,PerugiaMS4_HMP16b}.
\par

Let $A\in\mathbb{C}^{N_h\times N_h}$ denote the matrix associated with the sesquilinear form $\mathcal{A}_h(\cdot,\cdot)$,
and $\vec{b}\in\mathbb{C}^{N_h}$ the vector associated with the functional $\ell_h(\cdot)$, for $N_h:=\mbox{dim}(\mathrm{PW}_p(\mathcal{T}_h))$.
Then, the algebraic linear system associated with the PWDG method \eqref{perugia_pwdg} on the mesh $\mathcal{T}_h$ is
$A \vec{u} = \vec{b}$.

\section{Conditioning of the plane wave basis}
In this section, we investigate numerically the conditioning
of the local plane wave basis.
In order to do so, we consider the spectral condition number
of the local mass matrix $M_K\in\mathbb{C}^{p\times p}$ on a single element $K\in\mathcal{T}_h$. 
From \cite{PerugiaMS4_Gittelson} we get $M_{K,jj}=|K|$, and
\begin{align*}
 M_{K,jl} &= \int_K e^{ik\vec{d}_j\cdot (\vec{x}-\vec{x}_K)}\overline{e^{ik\vec{d}_l\cdot (\vec{x}-\vec{x}_K)}} \, \mathrm{d}\vec{x} \\
 &= -\sum_{F\in\partial K\cap \mathcal{F}_h} \frac{ik(\vec{d}_j-\vec{d}_l)\cdot\vec{n}}{k^2(\vec{d}_j-\vec{d}_l)\cdot(\vec{d}_j-\vec{d}_l)}\int_Fe^{ik(\vec{d}_j-\vec{d}_l)\cdot(\vec{x}-\vec{x}_K)}\, \mathrm{d}s,
\end{align*}
for $j\neq l$, $1\leq j,l\leq p$, which can be evaluated in closed form.
Note that the entries of $M_K$ tend to $|K|$ as $k(\vec{d}_j-\vec{d}_l)\cdot(\vec{x}-\vec{x}_K)$ tends to zero;
hence, small values of the element size $h$ and wavenumber $k$, or a small angle between two plane wave directions,
lead to ill conditioning.

\subsection{Dependence on the shape}
\begin{figure}[tb]
\centering
\includegraphics[width=0.49\textwidth]{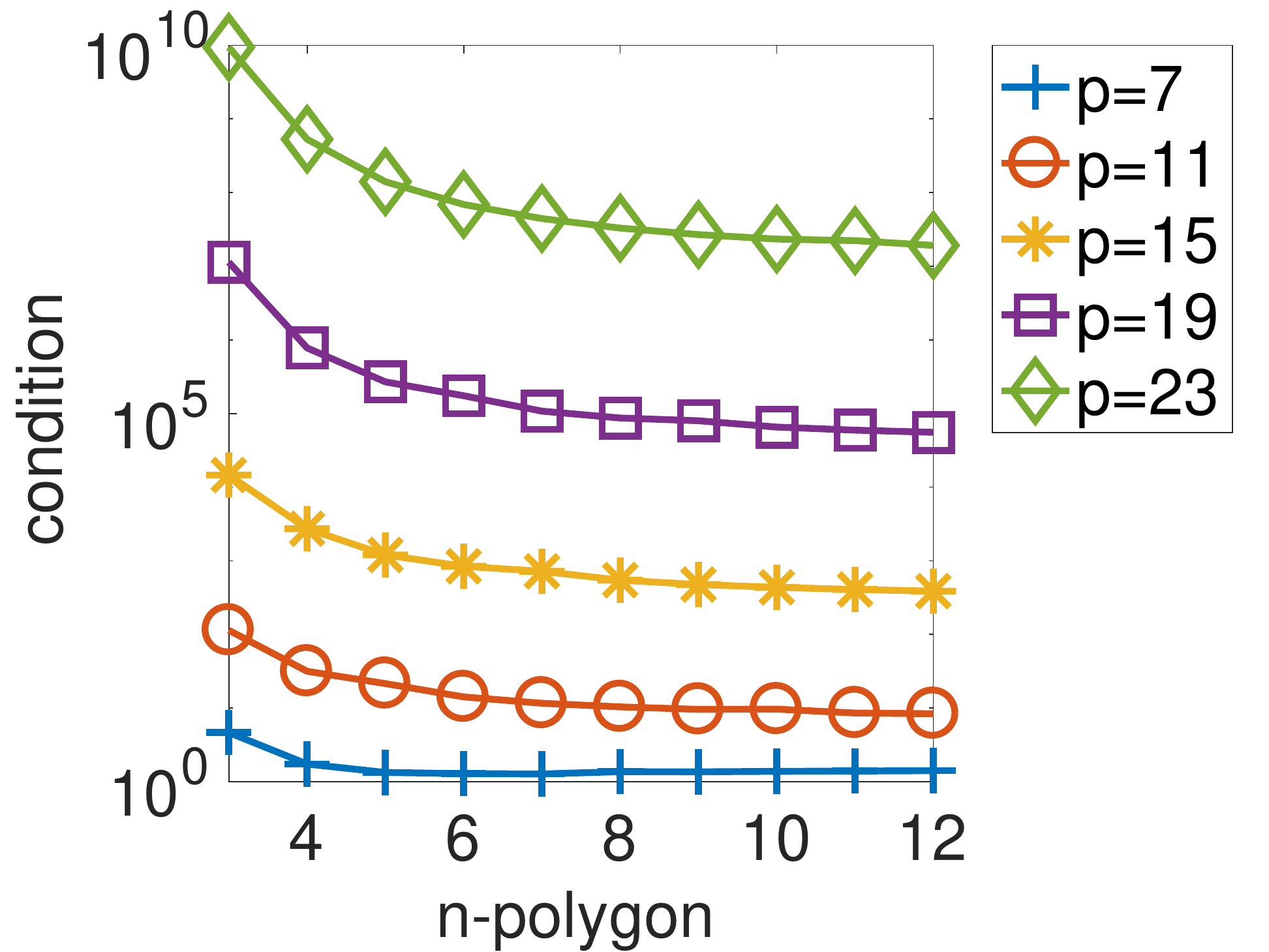}
\includegraphics[width=0.49\textwidth]{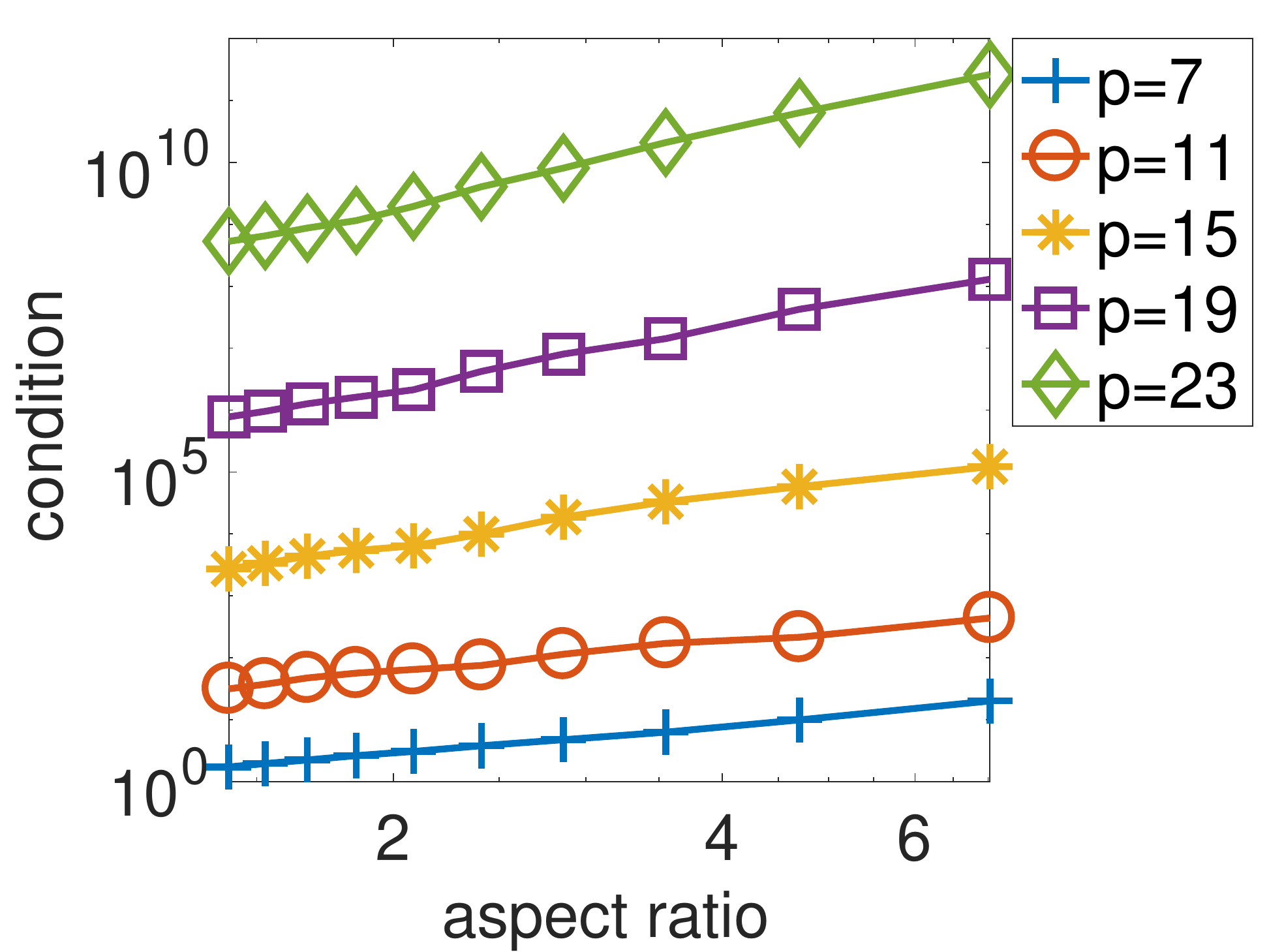}
\caption{Spectral condition numbers of $M_K$ for regular $n$-polygons (left) and anisotropic rectangles (right) with $h=1$ and $k=10$.}
\label{perugia_fig:1}
\end{figure}
The initial numerical experiments investigate the conditioning of the basis
for different shapes of the element,
for a fixed wavenumber $k=10$.

Firstly, we consider regular $n$-polygons with element size $h=1$; cf. Figure~\ref{perugia_fig:1} (left).
We observe that the condition numbers grow in 
the number of plane waves directions $p$, but are smaller for larger $n$.
In particular, the condition number is decreasing in the number
of sides $n$ and is asymptotically stable; hence, small
edges do not cause any problems.
It has been noted in \cite{PerugiaMS4_LHM13} that the conditioning of the basis depends on the
aspect ratio of the elements. Therefore, we consider 
a single anisotropic rectangle,
with size $h=1$, and vary its aspect ratio. 
We see that the condition number increases exponentially as the aspect ratio increases; cf. Figure~\ref{perugia_fig:1} (right).

\subsection{Dependence on $hk$ and $p$}
\begin{figure}[b]
\centering
\includegraphics[width=0.327\textwidth]{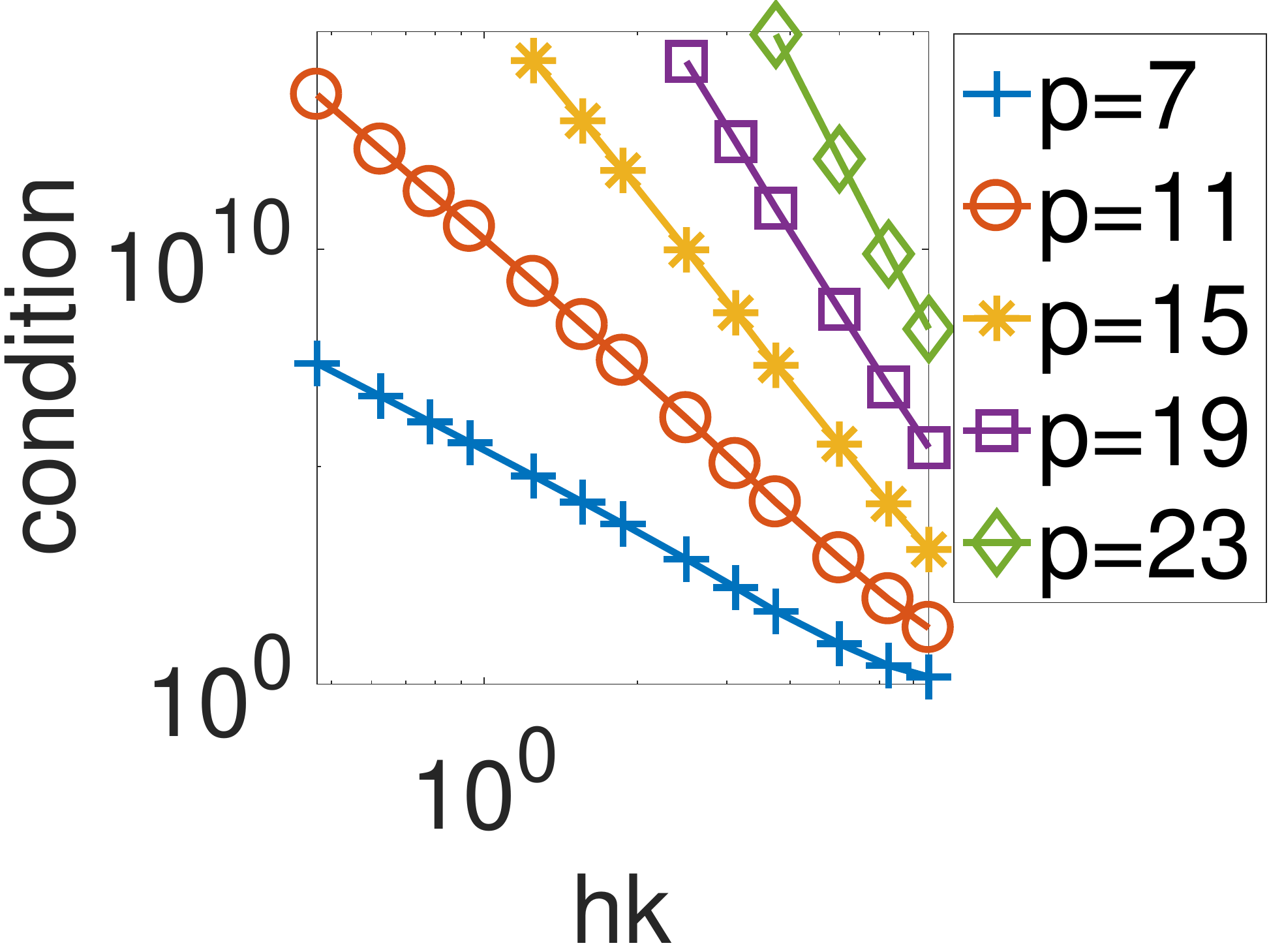}
\includegraphics[width=0.327\textwidth]{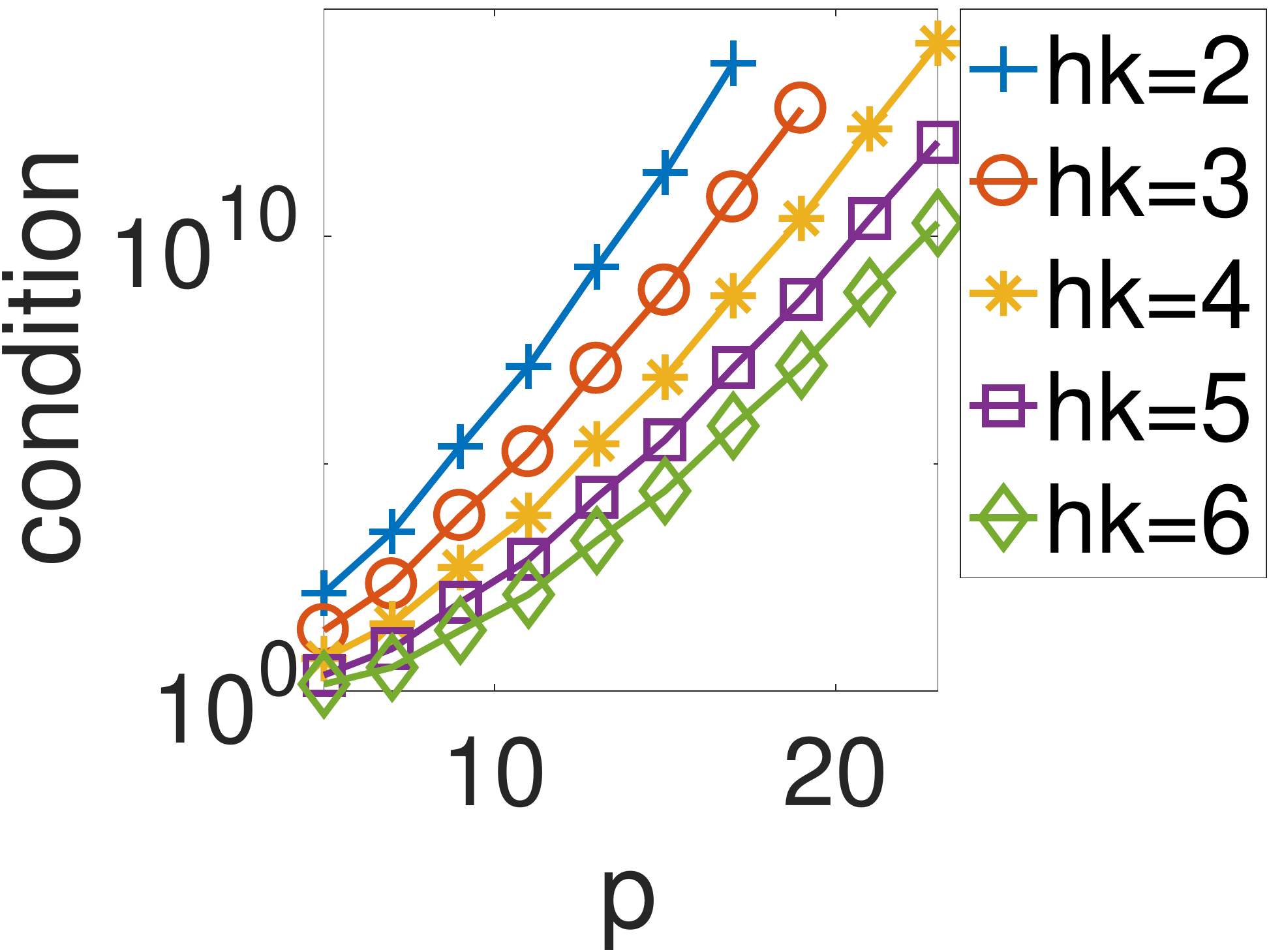}
\includegraphics[width=0.327\textwidth]{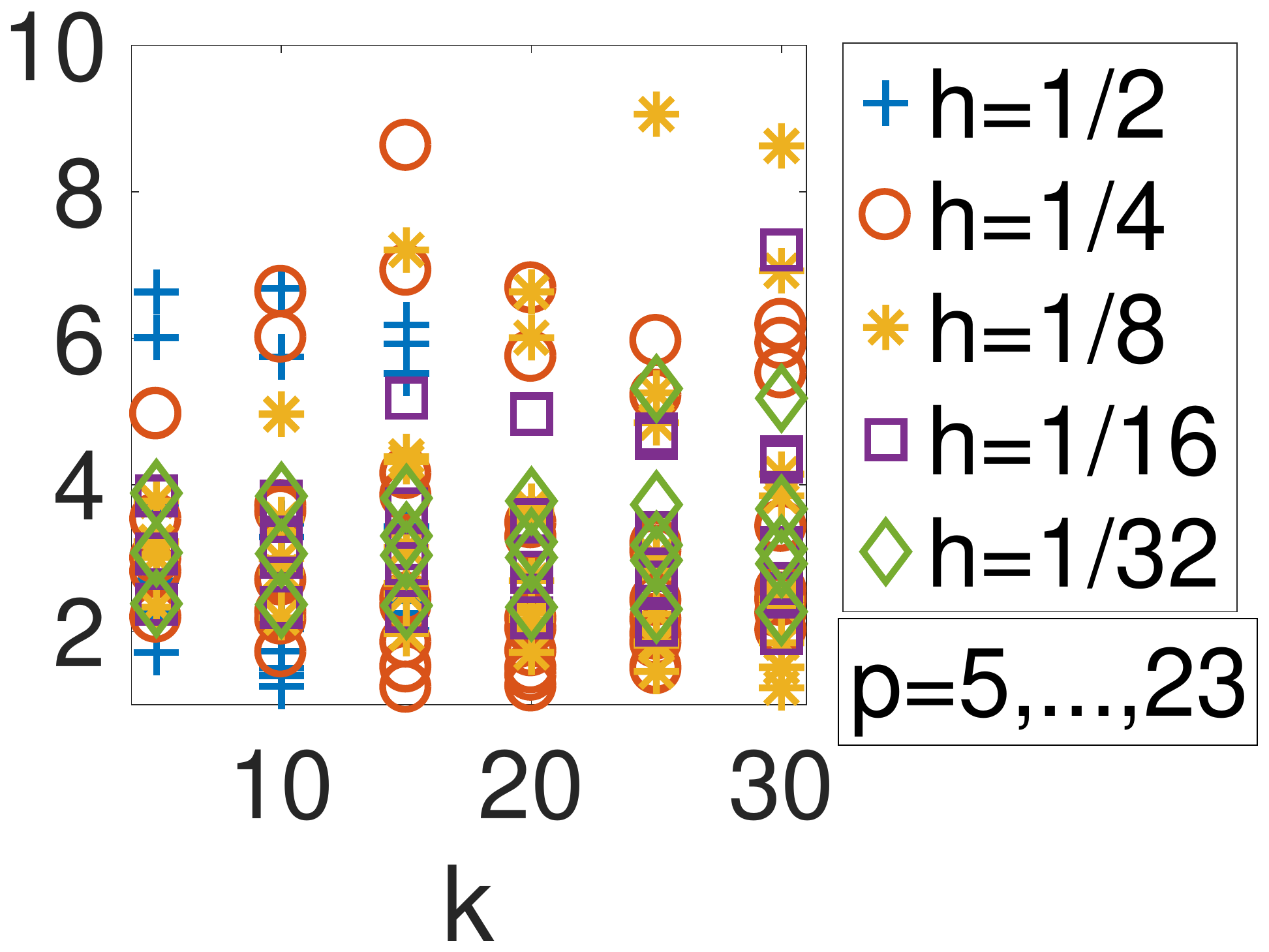}
\caption{Dependence of the condition number of $M_k$ on $hk$ and $p$, and verification of the approximation \eqref{perugia_cond}.}
\label{perugia_fig:2}
\end{figure}
In this section we empirically determine the dependence
of the condition number on $hk$ and $p$.
We restrict to the case of a single square element.
The numerical experiments displayed in the first two 
graphs in Figure~\ref{perugia_fig:2} suggest that the condition number is algebraic with respect to $hk$
and exponential with respect to $p$.
To get a more precise answer, we fitted the data obtained from numerous numerical experiments to
\begin{align}\label{perugia_cond}
\mbox{cond}_2(M_K)\approx \frac{2.34^{p\ln p}}{(hk)^{p-1}}.
\end{align}
In Figure~\ref{perugia_fig:2} (right) we show the values of 
$\mbox{cond}_2(M_K)/(\frac{2.34^{p\ln p}}{(hk)^{p-1}})$
for values $h=2^{-1},...,2^{-5}$, $k=5,\ldots,30$ and $p=5,\ldots,23$.
To obtain reliable data, we only plot data points for which $\mbox{cond}_2(M_K) < 10^{15}$,
due to double precision limitations, and for which $hk<10$, due to the resolution condition. 
Hence, we could only cover a moderate range of values
for $h$, $k$ and in particular $p$.
All the presented values of 
$\mbox{cond}_2(M_K)/(\frac{2.34^{p\ln p}}{(hk)^{p-1}})$
are between $1$ and $10$
(recall that the corresponding values of $\mbox{cond}_2(M_K)$ are between $1$ and $10^{15}$),
which confirms that the approximation \eqref{perugia_cond} is
reasonable, at least
for moderate $h$, $k$ and $p$.

\section{Orthogonalization of the plane wave basis}
\begin{figure}[b]
\centering
\includegraphics[width=0.32\textwidth]{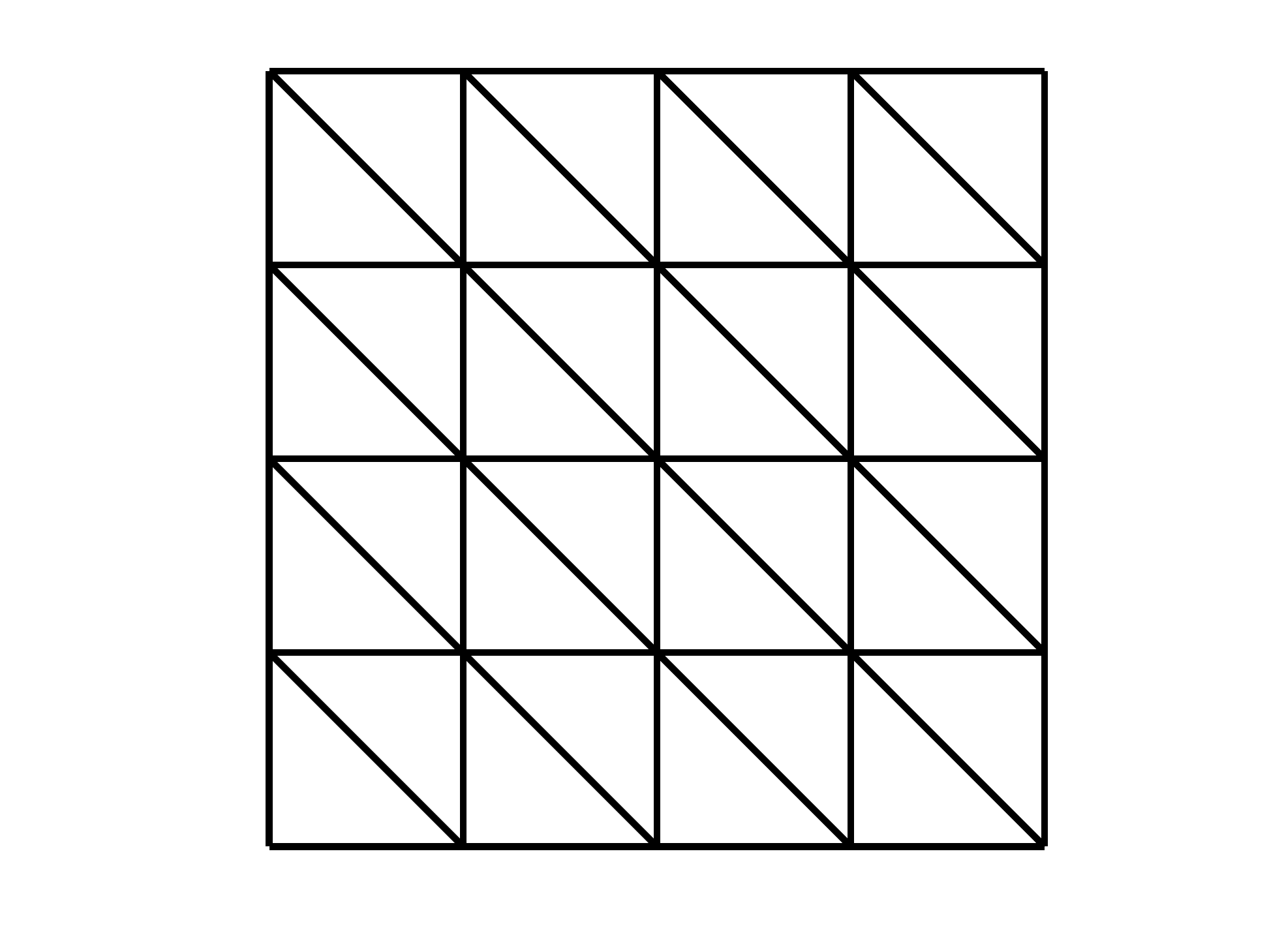}
\includegraphics[width=0.32\textwidth]{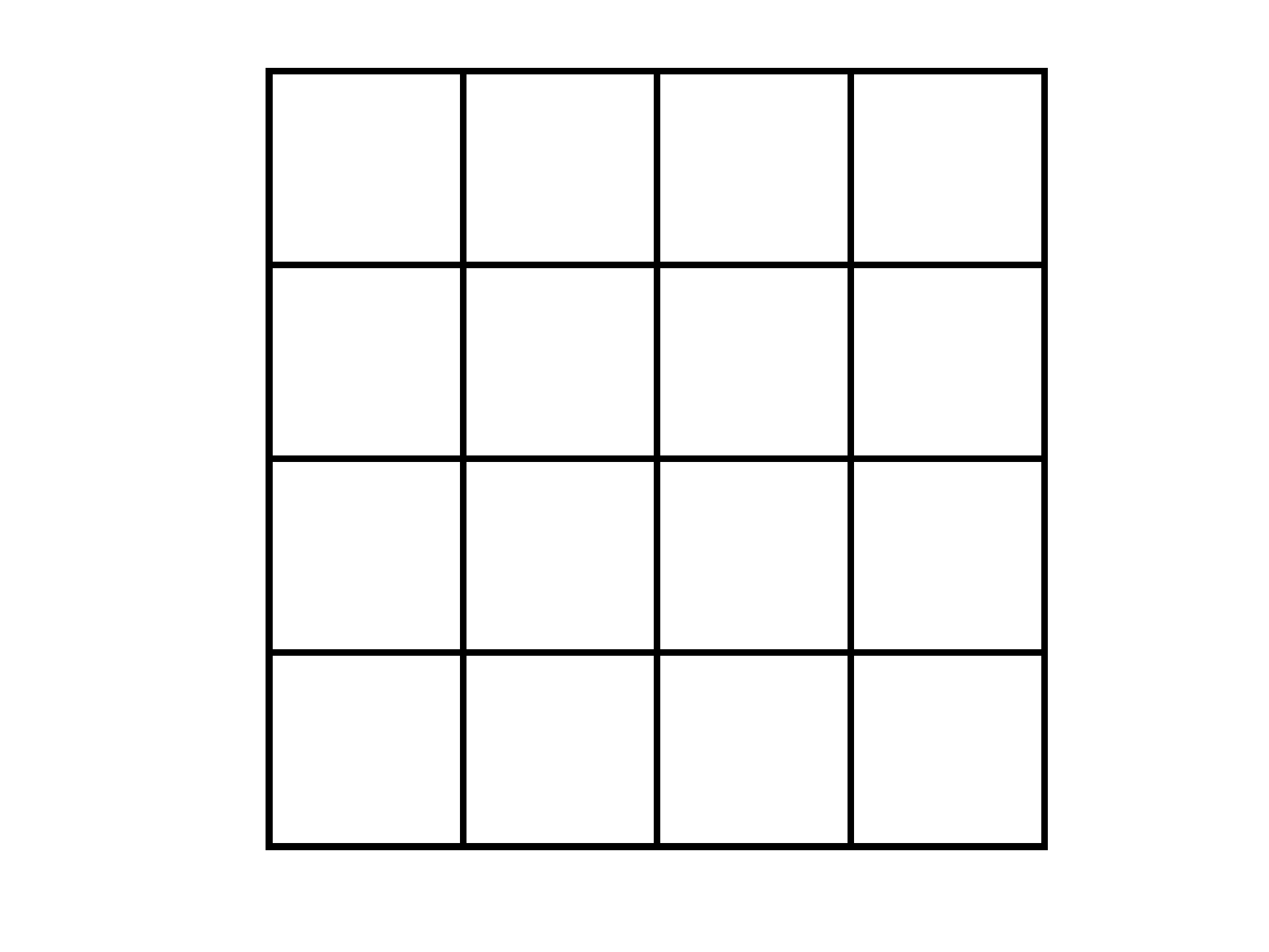}
\includegraphics[width=0.32\textwidth]{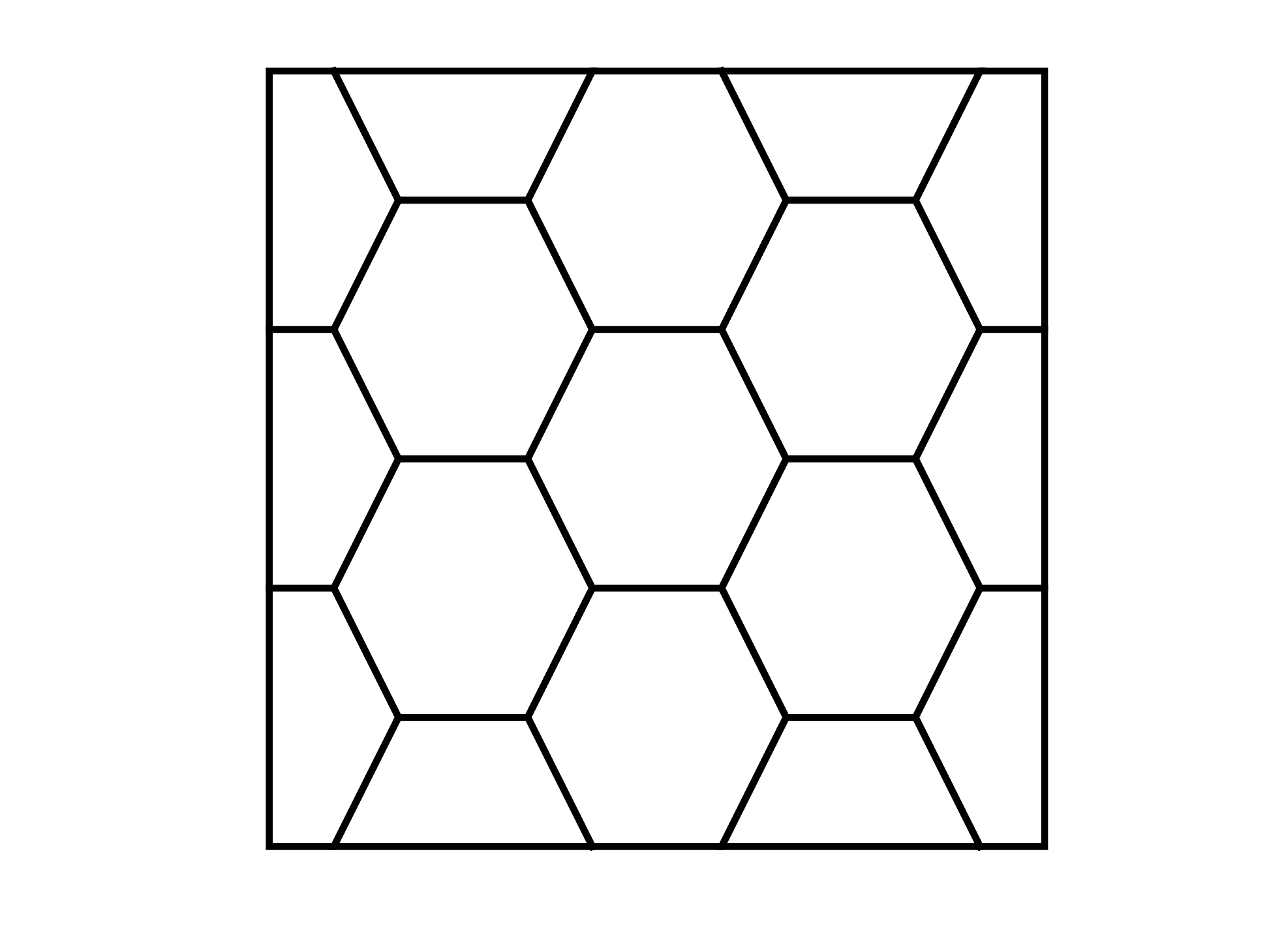}
\caption{Three different meshes.}
\label{perugia_meshes}
\end{figure}
In the previous section, we have observed that the
condition number of the local basis is large
for small $hk$ or large $p$.
In this section we aim at improving the conditioning of the local basis
in order to lower the condition number of the global system matrix $A$.
Therefore, we will investigate the effect of orthogonalization of
the (local) basis functions on the conditioning of the (global) system matrix $A$.
A different approach has been presented in
\cite{PerugiaMS4_HMK02}, where 
improvement of the conditioning of the global system
is achieved by suitably designed
non-uniform distributions of $p$.
\par

\begin{figure}[tb]
\centering
\includegraphics[width=0.49\textwidth]{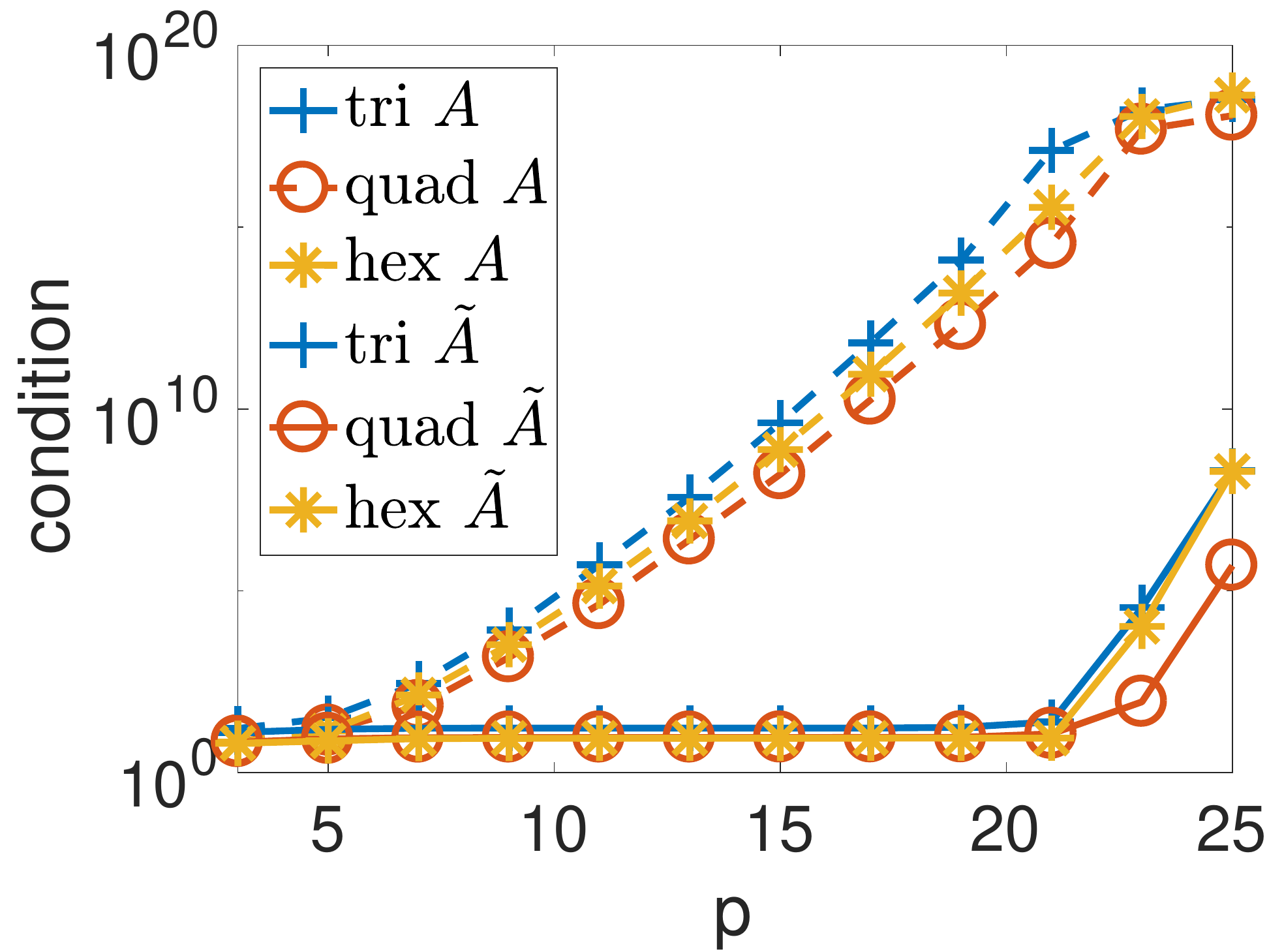}
\includegraphics[width=0.49\textwidth]{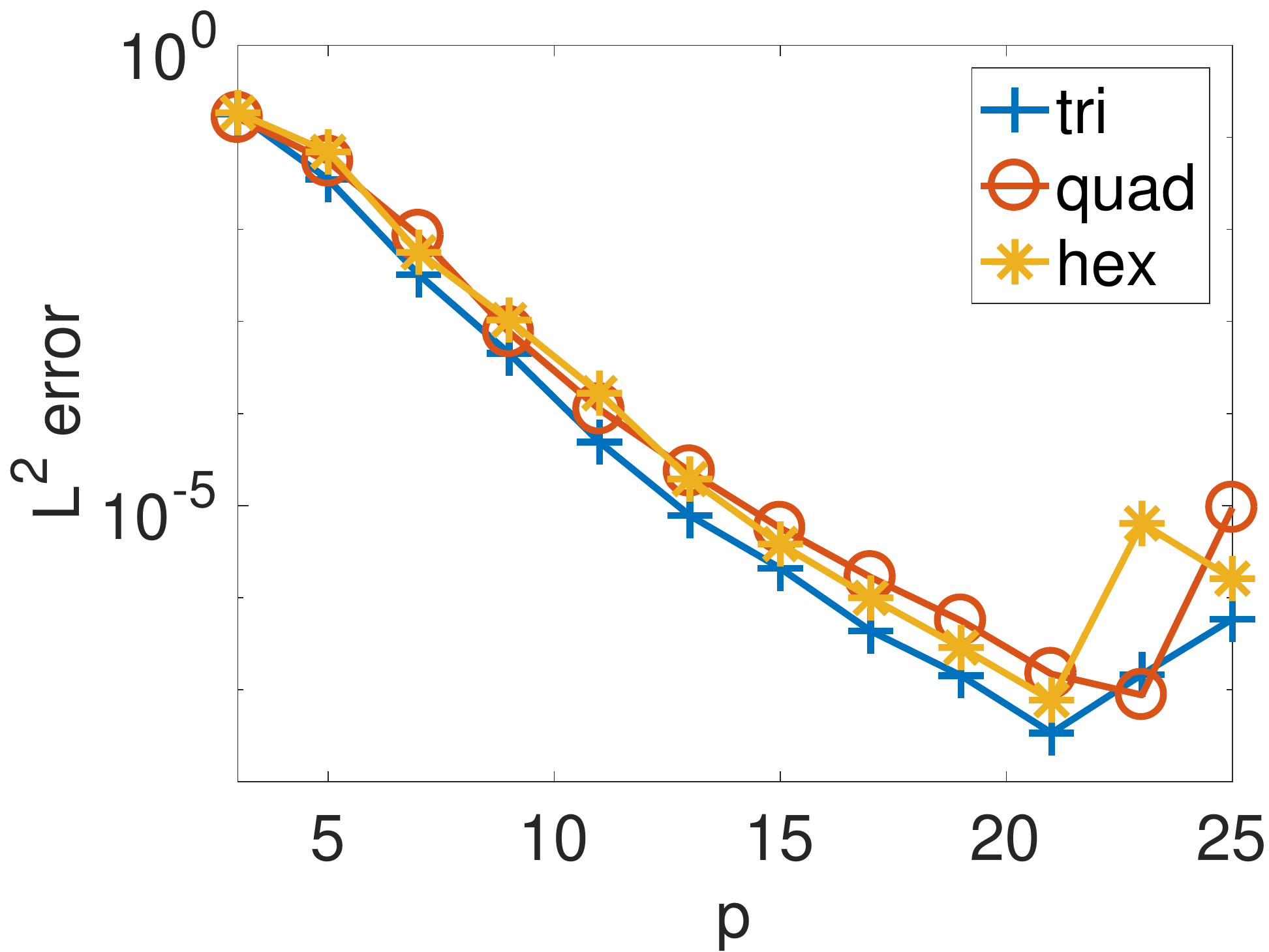}
\caption{Spectral condition numbers for $k=10$ in double precision arithmetic.}
\label{perugia_double}
\end{figure}
\begin{figure}[tb]
\centering
\includegraphics[width=0.49\textwidth]{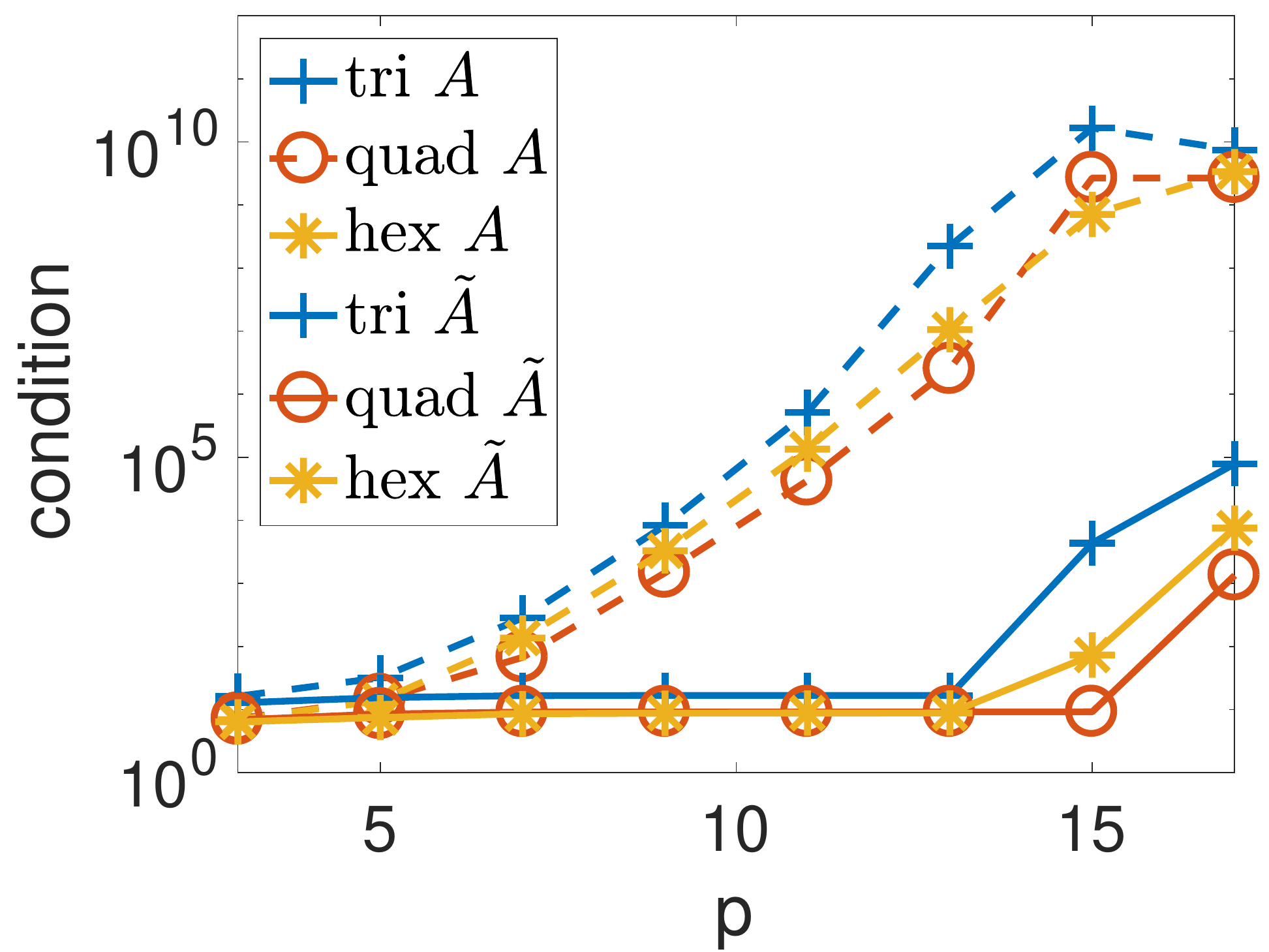}
\includegraphics[width=0.49\textwidth]{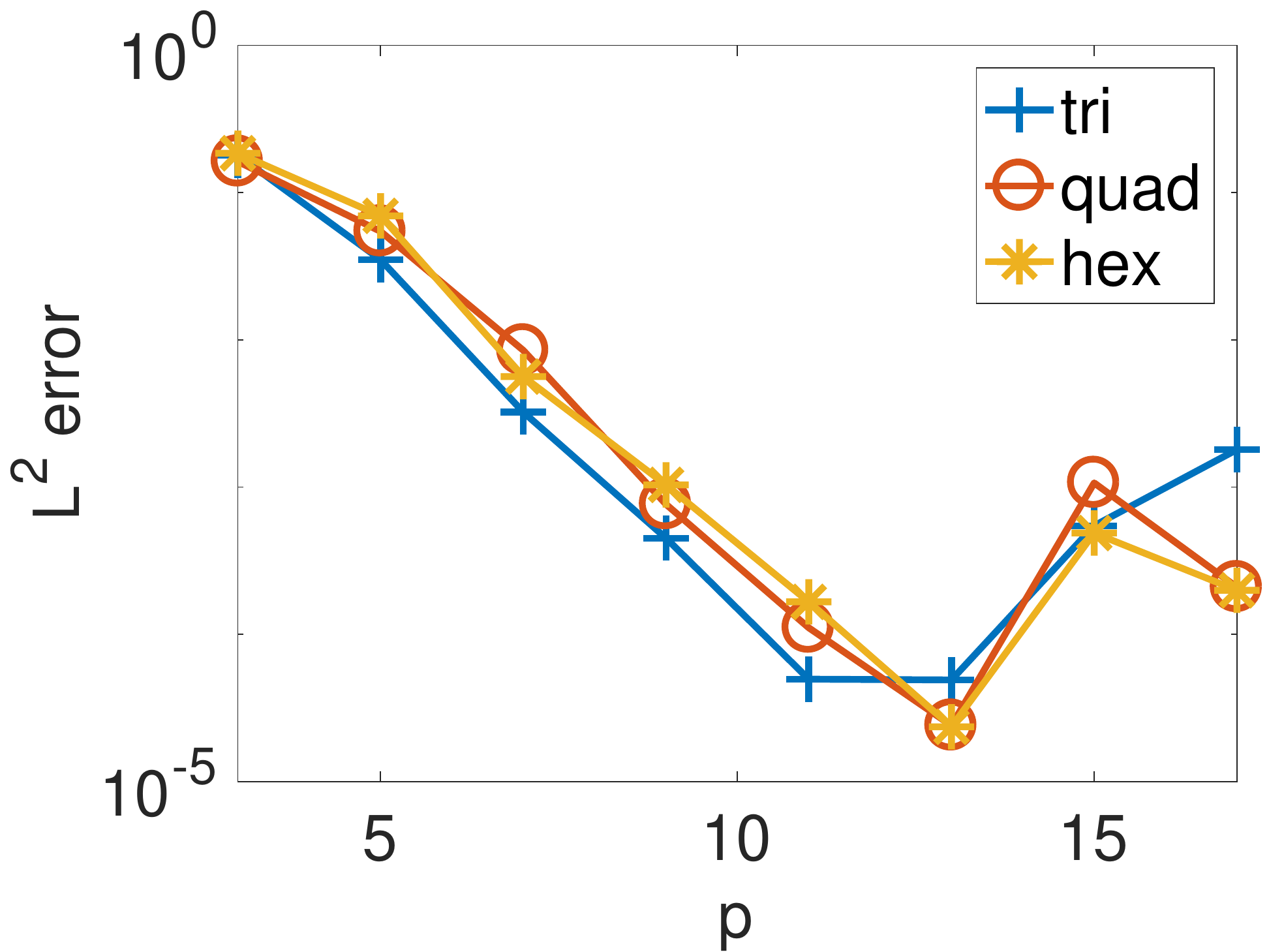}
\caption{Spectral condition numbers for $k=10$ in single precision arithmetic.}
\label{perugia_single}
\end{figure}
We compare the condition numbers of 
the system matrix $A$ with original basis functions with that of
the system matrix $\widetilde{A}:=Q^TAQ$ with orthogonalized basis functions
for the three different meshes displayed in 
Figure~\ref{perugia_meshes}. 
Here, $Q\in\mathbb{C}^{N_h\times N_h}$ denotes the transformation matrix obtained by
modified Gram-Schmidt orthogonalization \cite{PerugiaMS4_Steward} of the (local) basis functions
with respect to the Hermitian part of the local system matrix $H(A_K):=\frac{A_K+\bar{A}_k^T}{2}$
on each element $K$ separately;
cf. \cite{PerugiaMS4_Bassi,PerugiaMS4_Mascotto,PerugiaMS4_Schweitzer}
for application of modified Gram-Schmidt to
partition of unity methods, (polynomial) DG methods, and 
virtual element methods, respectively.
In all experiments, we choose $k=10$ and only investigate the effect on the critical
dependence of the condition number on $p$.
As a model problem, we consider problem \eqref{perugia_problem} with 
$\Omega=(0,1)^2$, and
exact solution given by the Bessel function of the third kind (Hankel function)
$u(\vec{x}) = H^1_0(k\sqrt{(x_1+1/4)^2+x_2^2})$.
\par

In Figure~\ref{perugia_double} (left), we observe, for all meshes, the expected increase of the condition number in $p$
for the original system matrix $A$ (dashed lines), which results in the loss of accuracy in the $L^2$ error for $p>21$, cf.  Figure~\ref{perugia_double} (right),
when using a direct linear solver.
We observe major improvements of the condition numbers for the matrix $\widetilde{A}$ (solid lines) until $p=21$ 
when the (modified) Gram-Schmidt orthogonalization breaks down,
which directly correlates to
the point when the direct solver fails to produce a more accurate solution.
Note that, for the original matrix $A$, there is no such correlation.
To further investigate these results, we also carried out the same experiments
in single precision arithmetic.
Figure~\ref{perugia_single} shows the results for single precision, where
we observe that the loss 
of accuracy already occurs at $p=13$.
Note, again, that this loss is correlated to the failure of the orthogonalization,
indicated by the sudden increase of the condition numbers for $\widetilde{A}$ at $p=13$.

\begin{table}[t]
\small
\setlength{\tabcolsep}{0.8mm}
\begin{tabular}{|l|r|r|r|r|r|r|}\hline
                                              &  \multicolumn1{c|}{$p=5$} & \multicolumn1{c|}{$p=7$} & \multicolumn1{c|}{$p=9$} & \multicolumn1{c|}{$p=11$} & \multicolumn1{c|}{$p=13$} & \multicolumn1{c|}{$p=15$}  \\ \hline
$\lambda_{min}(H(A))$          & $7.75\cdot10^{-1}$ & $2.56\cdot10^{-1}$ & $1.16\cdot10^{-2}$ & $4.50\cdot10^{-4}$ & $6.73\cdot10^{-6}$ & $1.18\cdot10^{-7}$\\ 
$\lambda_{min}(H(A^{-1}))$   & $3.83\cdot10^{-2}$ & $2.70\cdot10^{-2}$ & $2.09\cdot10^{-2}$ & $1.71\cdot10^{-2}$ & $1.44\cdot10^{-2}$ & $1.25\cdot10^{-2}$\\  
GMRES ($A$)         &     60   &      96  &   135   &   159   &     193   &    217        \\     \hline
$\lambda_{min}(H(\widetilde{A}))$          & $1.10\cdot10^{-1}$ & $8.45\cdot10^{-2}$ & $8.29\cdot10^{-2}$ & $8.27\cdot10^{-2}$ & $7.44\cdot10^{-2}$ & $6.08\cdot10^{-2}$\\ 
$\lambda_{min}(H(\widetilde{A}^{-1}))$   & $5.02\cdot10^{-1}$ & $5.00\cdot10^{-1}$ & $5.00\cdot10^{-1}$ & $5.00\cdot10^{-1}$ & $5.00\cdot10^{-1}$ & $5.00\cdot10^{-1}$\\ 
GMRES ($\widetilde{A}$)     &     47   &      52  &     58   &     62   &       68   &      73      \\       \hline  
\end{tabular}
\caption{Eigenvalue approximations and GMRES iteration count for the original basis and the orthogonalized basis 
using the second mesh of Figure~\ref{perugia_meshes} and $k=10$.}
\label{perugia_table}
\end{table}

Finally, we are interested in the 
effect of the orthogonalization on the convergence of iterative solvers such as GMRES.
From the convergence theory of GMRES \cite{EES83}, provided that $H(A)$,
the Hermitian part of the system matrix A, is positive definite, the residual
contraction factor, for the residual $r_j$ at iteration $j$, can be bounded as
\[
\frac{\|r_j\|}{\|r_0\|} \leq \left( 1 - \lambda_{min}(H(A))\lambda_{min}(H(A^{-1}))\right)^{j/2}.
\]
Therefore, in Table~\ref{perugia_table}, we report
$(\lambda_{min}(H(A)),\lambda_{min}(H(A^{-1})))$, 
for the system matrix $A$ with the original basis, and
$(\lambda_{min}(H(\widetilde{A})),\lambda_{min}(H(\widetilde{A}^{-1})))$,
for the system matrix $\widetilde{A}$ obtained
with the orthogonal basis, along with the number of GMRES iterations
needed in order to reduce the residual of a factor $10^{-10}$.
We observe that, for increasing $p$, the values $\lambda_{min}(H(A^{-1}))$ and
$\lambda_{min}(H(\widetilde{A}^{-1}))$ are fairly constant, 
while the values $\lambda_{min}(H(A))$ and
$\lambda_{min}(H(\widetilde{A}))$ decrease significantly.
However, the values of 
$\lambda_{min}(H(\widetilde{A}))$ decrease much 
more slowly than those of  
$\lambda_{min}(H(A))$,
which results in a far slower 
increase of GMRES iterations 
for the orthogonalized basis
than for the original one.
As we merely change the basis and do not precondition the system matrix,
we cannot expect a constant number of GMRES iterations.

\section{Conclusions}
We have provided empirical evidence that, in 2D, the condition number
of the plane wave basis is stable
for large edge counts in regular polygons, and 
grows as the aspect ratio of the anisotropy of the elements
increases. We also observed its
algebraic dependence 
on the product $hk$, and exponential dependence on~$p$.
It has been demonstrated that the condition number of the global system matrix
can be significantly lowered by a local modified Gram-Schmidt orthogonalization
with respect to the Hermitian part of the local system matrix; 
this results in faster convergence of the GMRES solver.
The improvement of the conditioning in 3D will be considered in future work.

\section*{Acknowledgments}  
The authors have been funded by the Austrian Science Fund (FWF) through the
project P 29197-N32. The third author has also been funded by the FWF through
the project F 65.

\ifx\undefined\bysame
\newcommand{\bysame}{\leavevmode\hbox to3em{\hrulefill}\,}
\fi

\end{document}